\documentclass[3p]{elsarticle}

\journal{Applied Mathematics Letters}

\usepackage[hidelinks]{hyperref}
\usepackage{amsmath}
\usepackage{amsfonts}
\usepackage{amssymb}
\usepackage{array}
\usepackage{subfigure}
\usepackage{color}

\graphicspath{{/}}
\renewcommand{\vec}[1]{\ensuremath \boldsymbol{#1}}	
\newcommand{\Pe}{\ensuremath \textnormal{Pe}}	
\newcommand{\Cn}{\ensuremath \textnormal{Cn}}

\title{Cahn-Hilliard on Surfaces: A Numerical Study\tnoteref{th1}}
\tnotetext[th1]{This work has been supported by the National Science Foundation through the 
	Division of Chemical, Bioengineering, Environmental, and Transport Systems Grant \#1253739.}
\author[ub]{P.~Gera}
\author[ub]{D.~Salac\corref{cor1}}
\cortext[cor1]{Corresponding author}
\address[ub]{Department of Mechanical and Aerospace Engineering, University at Buffalo,	Buffalo, New York 14260-4400}

\begin{document}
	\begin{abstract}
The Cahn-Hilliard system has been used to describe a wide number of phase separation processes, from 
co-polymer systems to lipid membranes. In this work the convergence properties of 
a closest-point based scheme is investigated. In place of solving the original
fourth-order system directly, two coupled second-order systems are solved. The system
is solved using an approximate Schur-decomposition as a preconditioner. The results
indicate that with a sufficiently high-order time discretization the method
only depends on the underlying spatial resolution.
	\end{abstract}

	
	\maketitle
\section{Introduction}

Phase field models are an important numerical technique for capturing the growth
of micro-structures~\cite{Chen2002} and can be used
to model various solidification processes such as that of a binary
alloy~\cite{beckermann1999modeling}. Instead of explicitly tracking the
solid-liquid interface, an order parameter is introduced that varies
continuously over a thin transition layer, which approximates a sharp solid
interface with a diffuse one.  Using this approach, boundary conditions at the
interface can be replaced by localized body-force terms, simplifying the
conservation equations that govern the order parameter.  Even though the term
``phase field" was coined much later, the basic idea has been in existence since
the 1800s when van der Waals studied fluid interfaces using a density
field~\cite{rowlinson1979translation}. Phase field models have also been used to
study many physical applications including defects~\cite{karma2001phase},
deformations~\cite{wang2010phase}, and the growth of
dendrites~\cite{warren1995prediction}. A variety of numerical
methods have been employed, including finite
differences~\cite{saylor2007diffuse},
finite elements~\cite{gomez2008isogeometric},
and spectral methods~\cite{chen1998applications}. 
One well known conservative method to capture phase dynamics is the Cahn-Hilliard
model. The original purpose of this model was to investigate the phase
transition of binary alloys~\cite{cahn1958free} and it has been used to
model the phase dynamics of polymers and ceramics~\cite{cogswell2010phase}, thermo-capillary
flows~\cite{verschueren2001diffuse}, 
the Rayleigh-Taylor instability~\cite{lowengrub1999topological},
droplet breakup~\cite{jacqmin1996energy}, liquid-liquid jets pinching
off~\cite{longmire1999comparison}, multicomponent lipid
vesicles~\cite{funkhouser2014dynamics}, and for the tracking of tumor
growth~\cite{wise2008three}.

Solving the Cahn-Hilliard equation on a curved surface is a difficult task for 
a number of reasons. It is a fourth order partial differential equation which
is numerically stiff, requiring small time steps for stability.
In this work, the fourth-order PDE is split into a set of coupled
second-order differential equations. The novelty of this work
lies in the use of the Closest-Point method~\cite{doi:10.1137/130929497}
and the associated numerical convergence analysis. This particular
formulation of the closest point method was originally applied to
second-order partial differential equations, while the splitting 
scheme here allows for the solution of fourth-order differential equations.
Additionally, a qualitative spatial analysis is also shown
to demonstrate that sufficient spatial resolution
is required to capture all aspects of the dynamics.

\section{Governing Equations}

Consider the boundary of a region, $\Gamma$, with an outward facing normal $\vec{n}$.
The non-dimensional surface Cahn-Hilliard equations describe the temporal evolution
of an order parameter $f(\vec{x},t)$ on the surface $\Gamma$, where the evolution is driven by
the gradient of the chemical potential,
\begin{align} 
	\frac{\partial f}{\partial t} = \dfrac{1}{\Pe} \nabla_s \cdot (\nu\nabla_s \mu),\label{eq:CH_continuity}
\end{align} 
where $\nabla_s$ is the surface gradient operator, $\nu$ is the mobility, and 
$\Pe$ is the surface Peclet number which relates the strength of any surface advection to diffusion.
The chemical potential, $\mu$, is the 
variational derivative of the
surface free energy,
\begin{align} 
	E_s=\int_\Gamma\left(g(f)+\dfrac{\Cn^2}{2}\|\nabla_s f\|^2\right)\;dA,
\end{align} 
where $g$ is often taken as a double well potential and 
$\Cn$ is the Cahn number, which indicates the strength of the energy penalty associated
with surface boundaries. 
This results in the following expression for the chemical potential,
\begin{equation} 
	\mu = \frac{dg}{df}-{\Cn}^2\Delta_s f, \label{eq:chempot}
\end{equation} 
where $\Delta_s$ is the surface Laplacian operator.
Using the chemical potential expression in Eq. (\ref{eq:CH_continuity})
results in a fourth-order evolution equation for the phase field.

\section{Numerical Methods}

Assume that the mobility is constant and equal to 1
and that the equations are discretized on a standard Cartesian grid using finite difference
approximations. Let the interface be described by a level set Jet scheme~\cite{seibold2012jet}.
The Cahn-Hilliard system can be written as a pair of coupled, second-order differential
equations~\cite{lowengrub2007surface,Li2012}, 
\begin{align}
	\frac{\partial f}{\partial t} - \frac{1}{\Pe}\Delta_s  \mu = 0 \quad \textnormal{and} \quad \mu+\Cn^2\Delta_s f =g'(f).
\end{align}
Discretizing this in time the system can be written as
\begin{equation}
	\begin{bmatrix}
		\vec{I} & \Cn^2\vec{L}_s\\
		-\frac{\Delta t}{\Pe\beta_0}\vec{L}_s & \vec{I}
	\end{bmatrix}
	\begin{bmatrix}
		\vec{\mu}^{n+1}\\
		\vec{f}^{n+1}
	\end{bmatrix}
	=
	\begin{bmatrix}
		g'(\hat{\vec{f}})\\
		-\frac{\beta_1}{\beta_0} \vec{f}^n -\frac{\beta_2}{\beta_0} \vec{f}^{n-1}
	\end{bmatrix},
	\label{eq:CH_block}
\end{equation}
where $\vec{f}^n$ and $\vec{f}^{n-1}$ are the solutions at times
$t^n$ and $t^{n-1}$, respectively, with $\Delta t=t^{n}-t^{n-1}$ a constant time step,
and $\vec{L}_s$ is the discrete approximation to the surface Laplacian, $\Delta_s$.
The approximation to the solution at time $t^{n+1}$ is given by $\hat{\vec{f}}$ and
$\vec{I}$ is the identity matrix.
A first-order backward-finite-difference (BDF1) scheme is given by
$\beta_0=1$, $\beta_1=-1$, and $\beta_2=0$ with $\hat{\vec{f}}=\vec{f}^n$.
A second-order backward-finite-difference (BDF2) scheme is given by
$\beta_0=3/2$, $\beta_1=-2$, and $\beta_2=1/2$ with $\hat{\vec{f}}=2\vec{f}^n-\vec{f}^{n-1}$~\cite{Fornberg1988}.

The discretization of the surface Laplacian is performed using 
a version of the Closest-Point Method introduced by Chen 
and Macdonald~\cite{doi:10.1137/130929497}.
The basic idea is to extend the solution away
from the interface such that it is constant in the normal direction. With this
extension, it is possible to write a surface differential equation as a standard
differential equation in the embedding space. In this scheme points
far from the interface do not need to be considered. Note that in situations
where a grid point might be near multiple interfaces extension errors may
be introduced, but for time-evolving interfaces, such as those undergoing pinch-off,
this error will be temporary.

Let $\vec{E}_1$ be a linear polynomial interpolation operator and $\vec{E}_3$ be a cubic polynomial interpolation operator.
For any point $\vec{x}$ not on the interface these operators return
the value of a function at the interface point closest to $\vec{x}$. For example, the operation $\vec{E}_3\vec{f}$ 
returns the value of $f$ at the point on the interface closest to $\vec{x}$ using the cubic interpolation function. Using this notation,
Eq. (\ref{eq:CH_block}) is re-written as
\begin{equation}
	\begin{bmatrix}
		\vec{I} & \Cn^2\left[\vec{E}_1\vec{L}+\alpha\left(\vec{E}_3-\vec{I}\right)\right]\\
		-\frac{\Delta t}{\Pe\beta_0}\left[\vec{E}_1\vec{L}+\alpha\left(\vec{E}_3-\vec{I}\right)\right] & \vec{I}
	\end{bmatrix}
	\begin{bmatrix}
		\vec{\mu}^{n+1}\\
		\vec{f}^{n+1}
	\end{bmatrix}
	=
	\begin{bmatrix}
		g'(\hat{\vec{f}})\\
		-\frac{\beta_1}{\beta_0} \vec{f}^n -\frac{\beta_2}{\beta_0} \vec{f}^{n-1}
	\end{bmatrix},
	\label{eq:CP_CH_block}
\end{equation}
with $\alpha=6/h^2$ where $h$ is the uniform grid spacing and $\vec{L}$ represents
the Cartesian finite difference approximation to the standard Laplacian, $\Delta$. The addition of $\alpha$ term ensures that the solutions
are constant in the normal direction. If this extension holds then surface operators can be replaced
with standard Cartesian operators. See Chen and Macdonald for complete
details~\cite{doi:10.1137/130929497}. 

The block system shown in Eq. (\ref{eq:CP_CH_block}) is solved using the preconditioned 
Flexible GMRES algorithm available in PETSc~\cite{petsc-web-page,petsc-user-ref,petsc-efficient},
where the preconditioner is an approximate Schur inverse.
Let $\vec{L}_E=\vec{E}_1\vec{L}+\alpha\left(\vec{E}_3-\vec{I}\right)$. The preconditioner is then
\begin{equation}
	\vec{P}=
	\begin{bmatrix}
		\vec{I} & - \Cn^2\vec{L}_E\\
		0 & \vec{I}
	\end{bmatrix}
	\begin{bmatrix}
		\vec{I}& 0\\
		0 & \vec{\hat{S}}^{-1}
	\end{bmatrix}
	\begin{bmatrix}
		\vec{I} & 0 \\
		\frac{\Delta t}{\Pe \beta_0}\vec{L}_E & \vec{I}
	\end{bmatrix},
\end{equation}
with a Schur complement of $\vec{\hat{S}}=\vec{I}+\frac{\Cn^2\Delta t_{max}}{\Pe \beta_0}\vec{L}_E\vec{L}_E$,
where $\Delta t_{max}=10^{-2}$ is the maximum time-step considered here.
The inverse of the Schur complement is calculated using an explicit LU-decomposition using MUMPS~\cite{Amestoy2001,Amestoy2006}.
Note that for a time step of $\Delta t=10^{-2}$, this results in an exact method and the GMRES algorithm converges 
within one iteration. For smaller time steps the Schur complement is no longer exact and GMRES algorithm converges after two or three iterations.
Additionally, for small Cahn numbers, say $\Cn\sim 3h$, where $h$ is the grid spacing, the numerical 
stiffness of the system is greatly reduced, which might remove the need for such a preconditioner.

\section{Results}
First consider the quantitative convergence of the method 
by comparing to a higher-order scheme. Results here are presented on a 
unit sphere in a domain spanning $[-1.25,1.25]^3$ 
with uniform grid spacing $h$ and fixed time step $\Delta t$. 
The
initial condition is taken to be $\cos\left(\cosh\left(5 x z\right)-10y\right)$, where $(x,y,z)$ is a 
point on the sphere. The mixing energy is a simple double-well potential, $g(f)=(f^4)/4-(f^2)/2$,
while the parameters are set to $\Cn=0.1$ and $\Pe=1.0$.
The reference solution is calculated via the spectral spherical operators in the Chebfun
Matlab package using a uniform $256^2$ mesh~\cite{Chebfun14} and LIRK4, which is a fourth-order accurate scheme, using a time
step of $10^{-4}$~\cite{Montanelli2017}.
All $l_\infty$-errors are computed at a time of $t=0.5$.

\begin{figure}
	\centering
	\subfigure[Error in Phase Field]{
		\includegraphics[width=7.0cm]{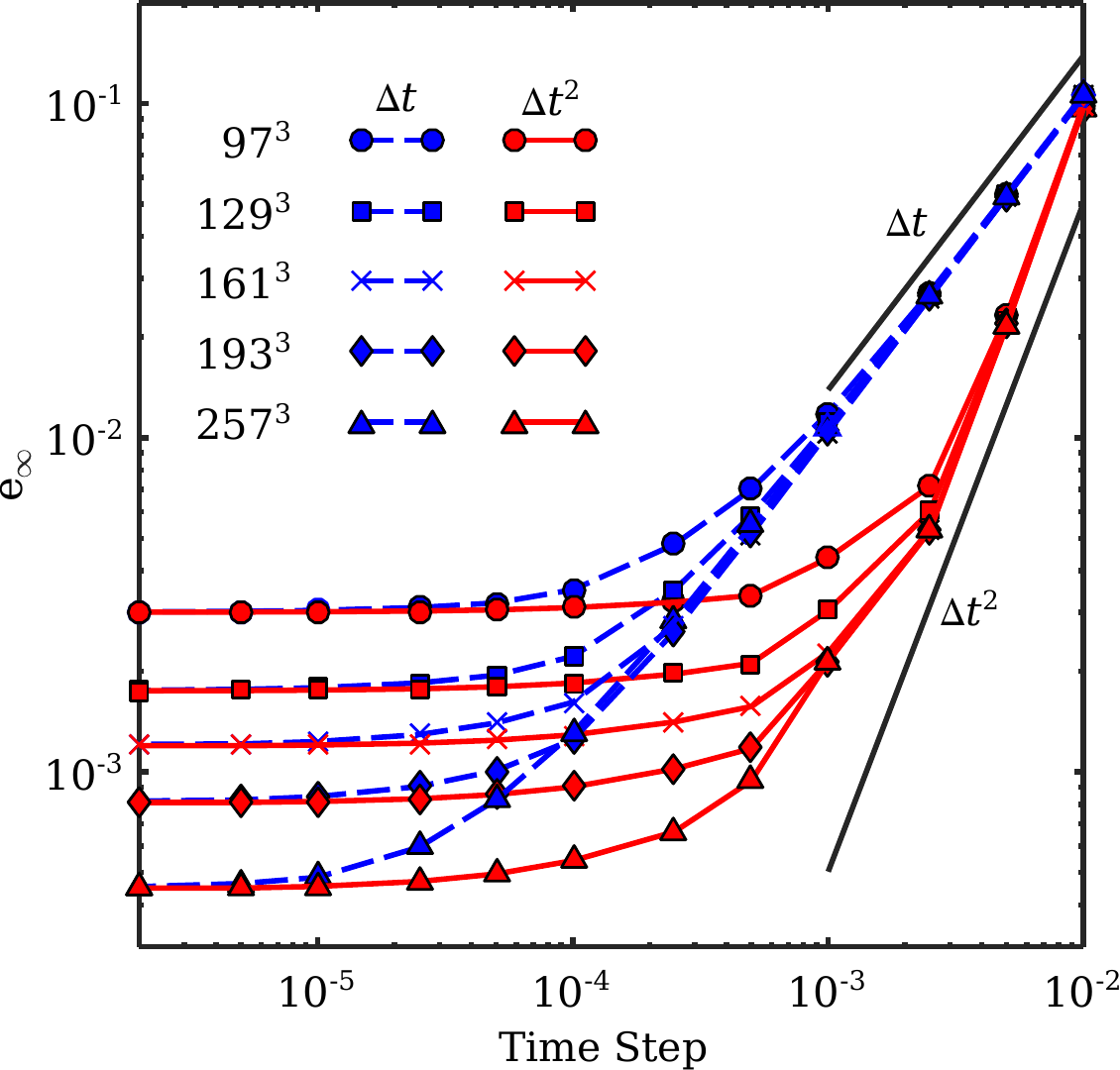}
    }
	\qquad
	\subfigure[Error in Chemical Potential]{
		\includegraphics[width=7.0cm]{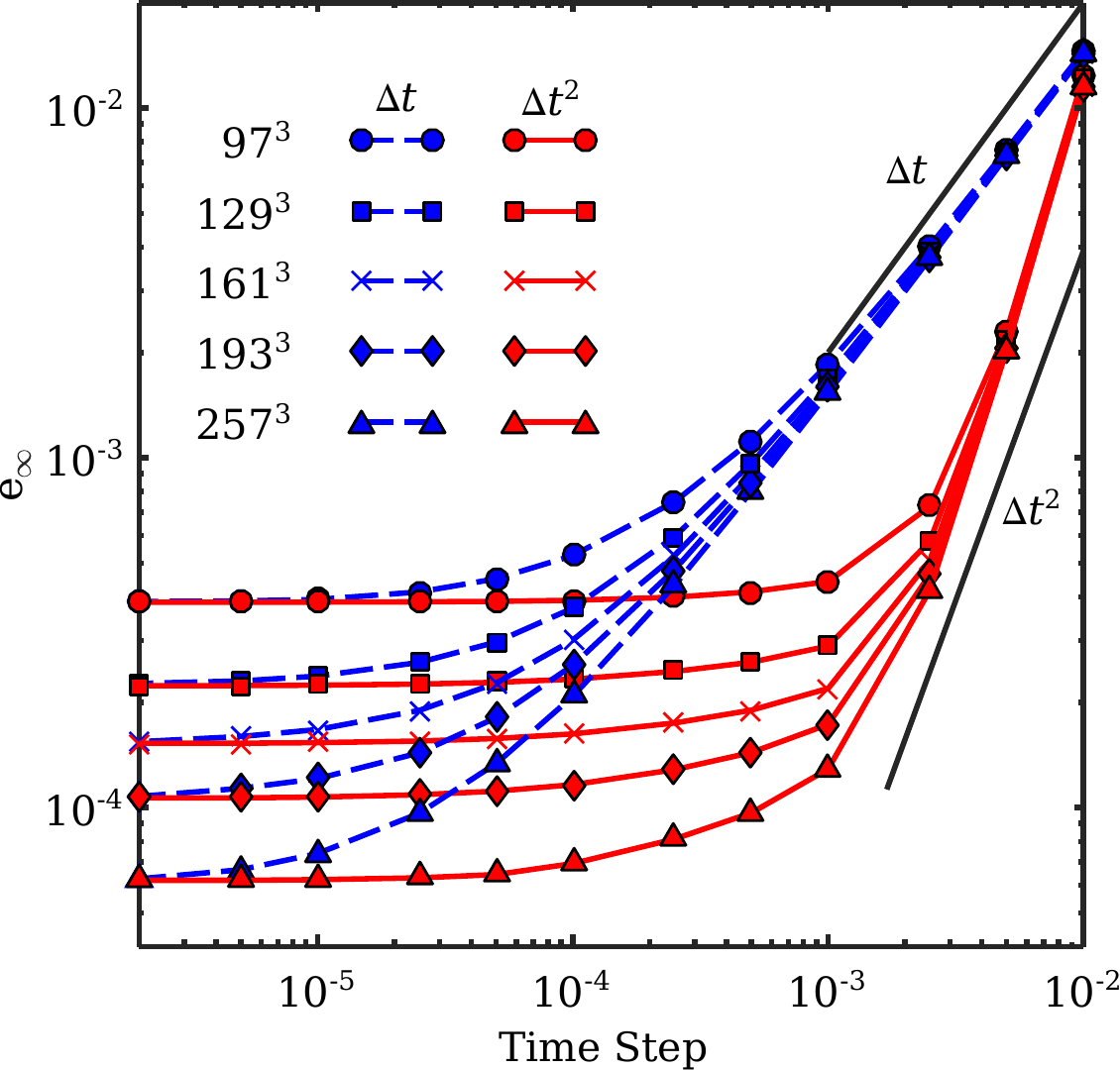}
	}	
	\caption{Time convergence results are presented for first and second order
    time schemes for different grid spacing. At lower time steps, error is
    independent of the time scheme used and is dominated by the grid size.}	
  \label{fig:timeConv}
\end{figure}

\begin{table}
	\centering
	\caption{Grid convergence results are presented using the methods of
	manufactured solutions using BDF2. Errors in the phase field and chemical potential
	at time of t=0.5 are shown using time step of $2.5 \times 10^{-6}$.}	
	\label{table:gridConv}
	
	\begin{tabular}{|c|c|c|c|c|c|}
		\cline{3-6} 
		\multicolumn{2}{c|}{} & \multicolumn{2}{c|}{Phase Field} & \multicolumn{2}{c|}{Chemical Potential}\tabularnewline
		\hline 
		N$^3$ & h & e$_\infty$ & Order & e$_\infty$ & Order \tabularnewline
		\hline 
		97 & 0.0260 & 3.0051$\times$10$^{-3}$ &      & 3.8621$\times$10$^{-4}$ & \tabularnewline
		\hline 
		129 & 0.0195 & 1.7484$\times$10$^{-3}$ & 1.88 & 2.2288$\times$10$^{-4}$ &1.91\tabularnewline
		\hline 
		161 & 0.0156 & 1.1978$\times$10$^{-3}$ & 1.70 & 1.5271$\times$10$^{-4}$ &1.69\tabularnewline
		\hline 
		193 & 0.0130 & 8.1018$\times$10$^{-4}$ & 2.14 & 1.0656$\times$10$^{-4}$ &1.97\tabularnewline
		\hline 
		257 & 0.0098 & 4.4920$\times$10$^{-4}$ & 2.05 & 6.1937$\times$10$^{-5}$ &1.89\tabularnewline		
		\hline 
		
	\end{tabular}	
\end{table}

\begin{figure}
	\centering
	\begin{tabular}{>{\centering}m{1.0cm} >{\centering}m{3.0cm}>{\centering}m{3.0cm} >{\centering}m{3.0cm} >{\centering}m{3.0cm}}
		& \multicolumn{4}{c}{Grid Size} \tabularnewline
		\cline{2-5}
		Time & $97^3$ & $129^3$ & $161^3$ & $257^3$\tabularnewline
		0.2 & \includegraphics[width=2.75cm]{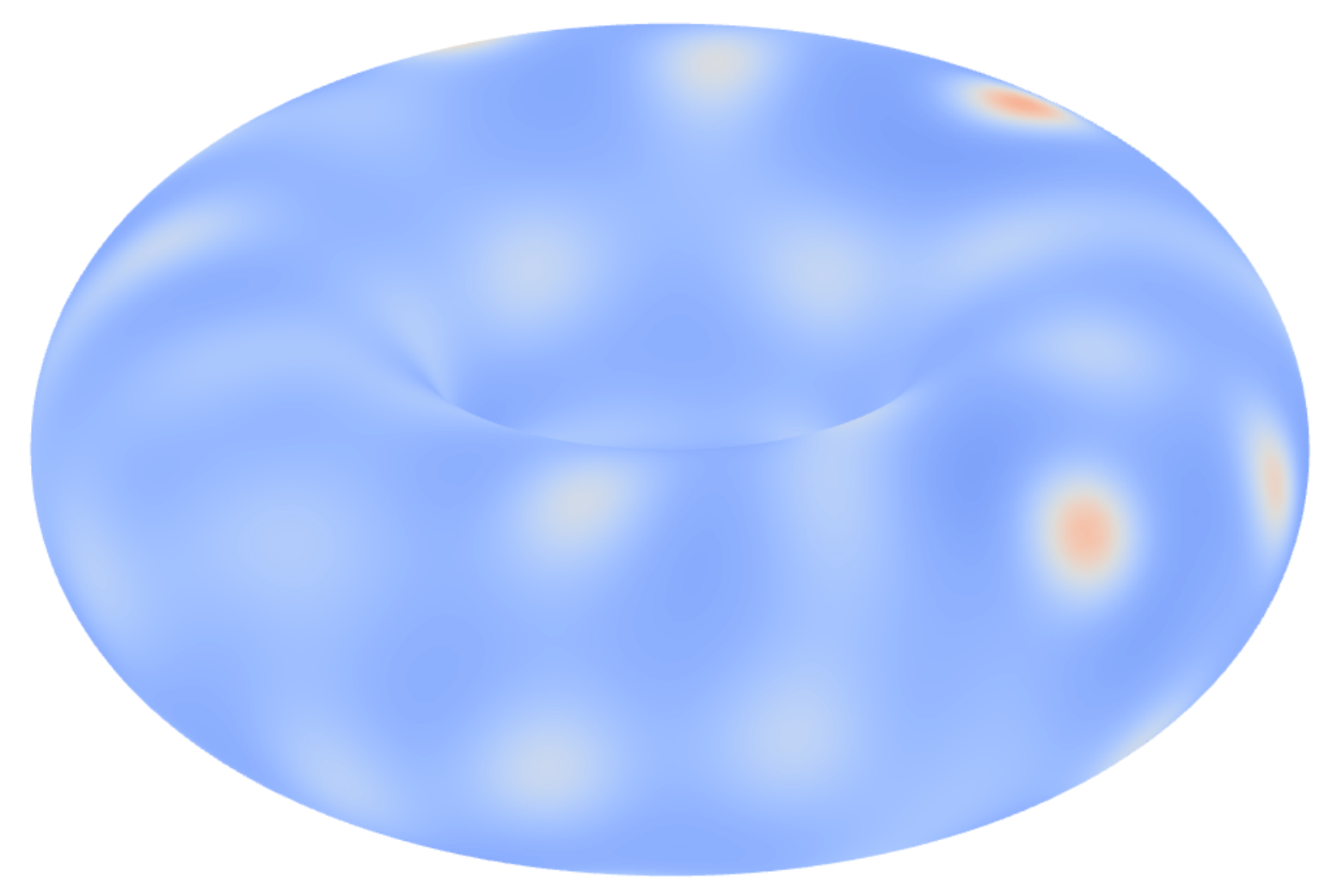} & \includegraphics[width=2.75cm]{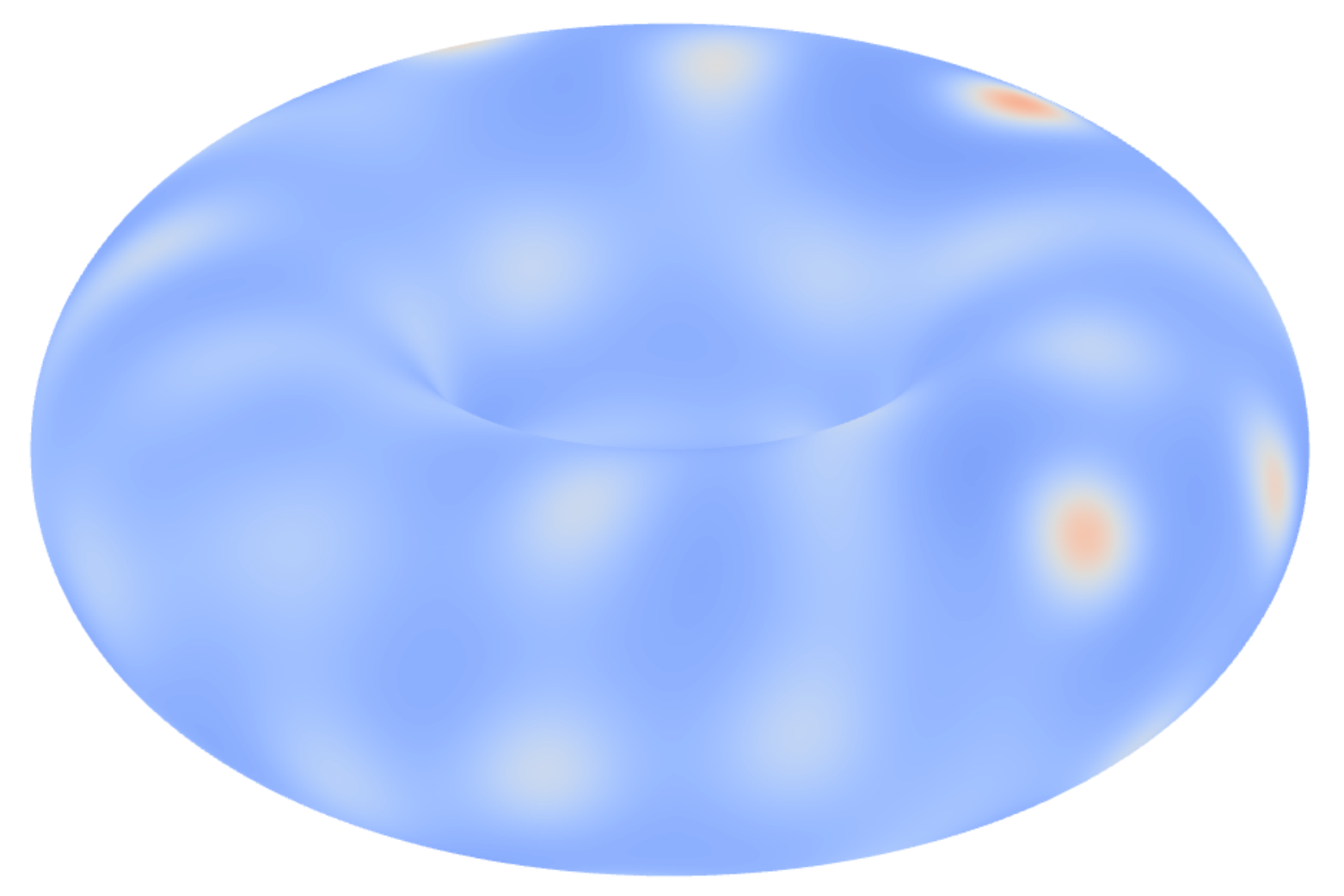} & \includegraphics[width=2.75cm]{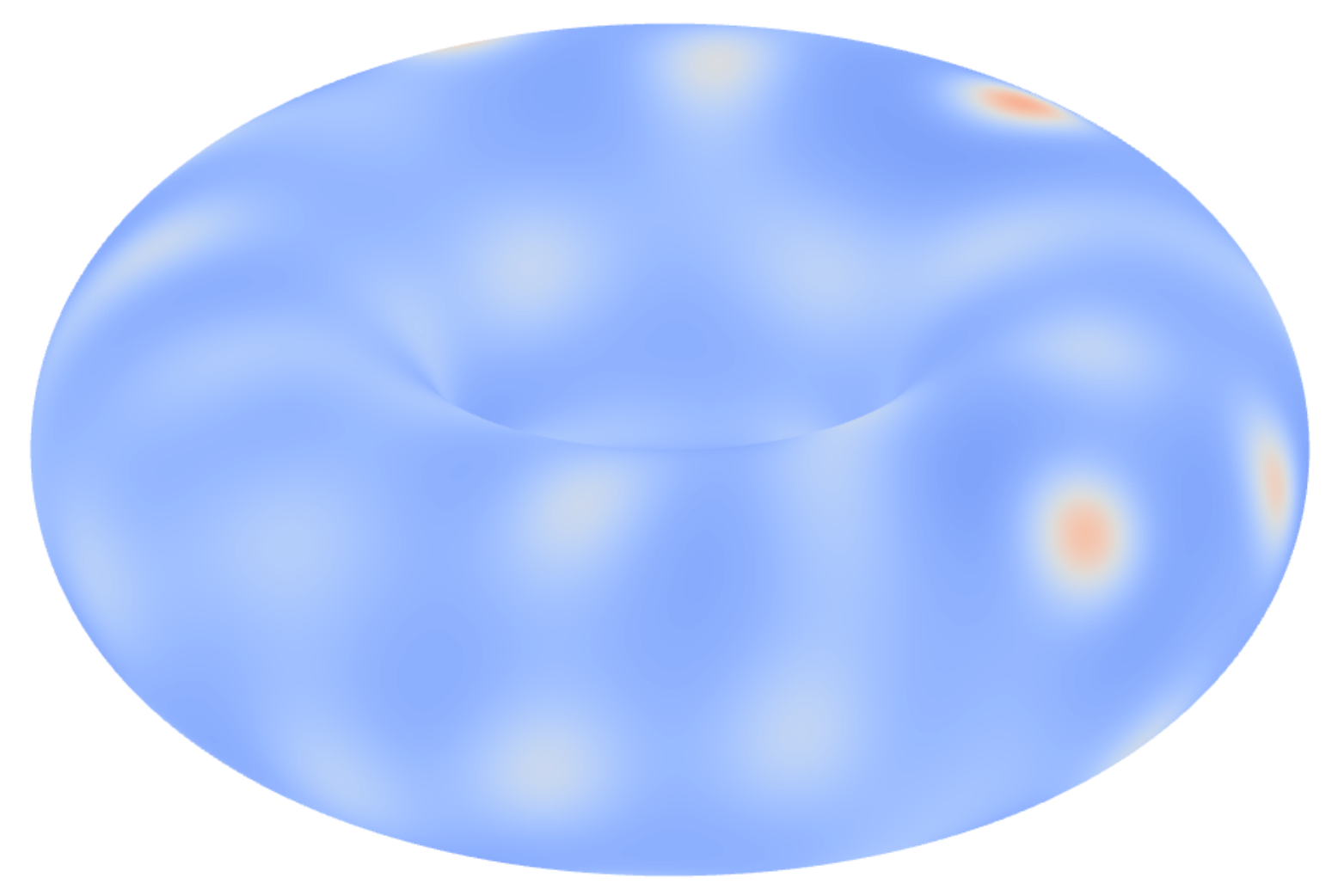} & \includegraphics[width=2.75cm]{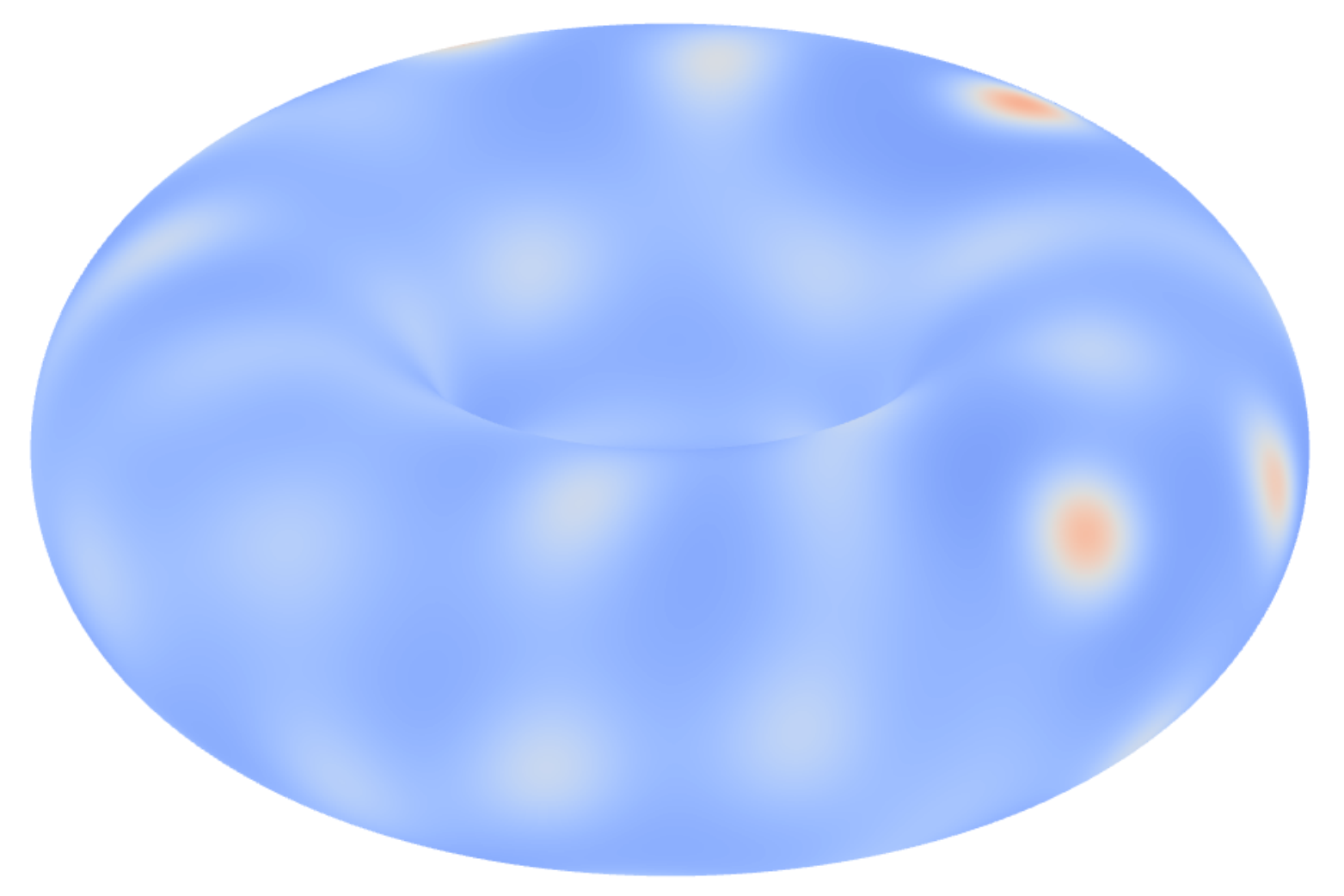} \tabularnewline
		0.3  & \includegraphics[width=2.75cm]{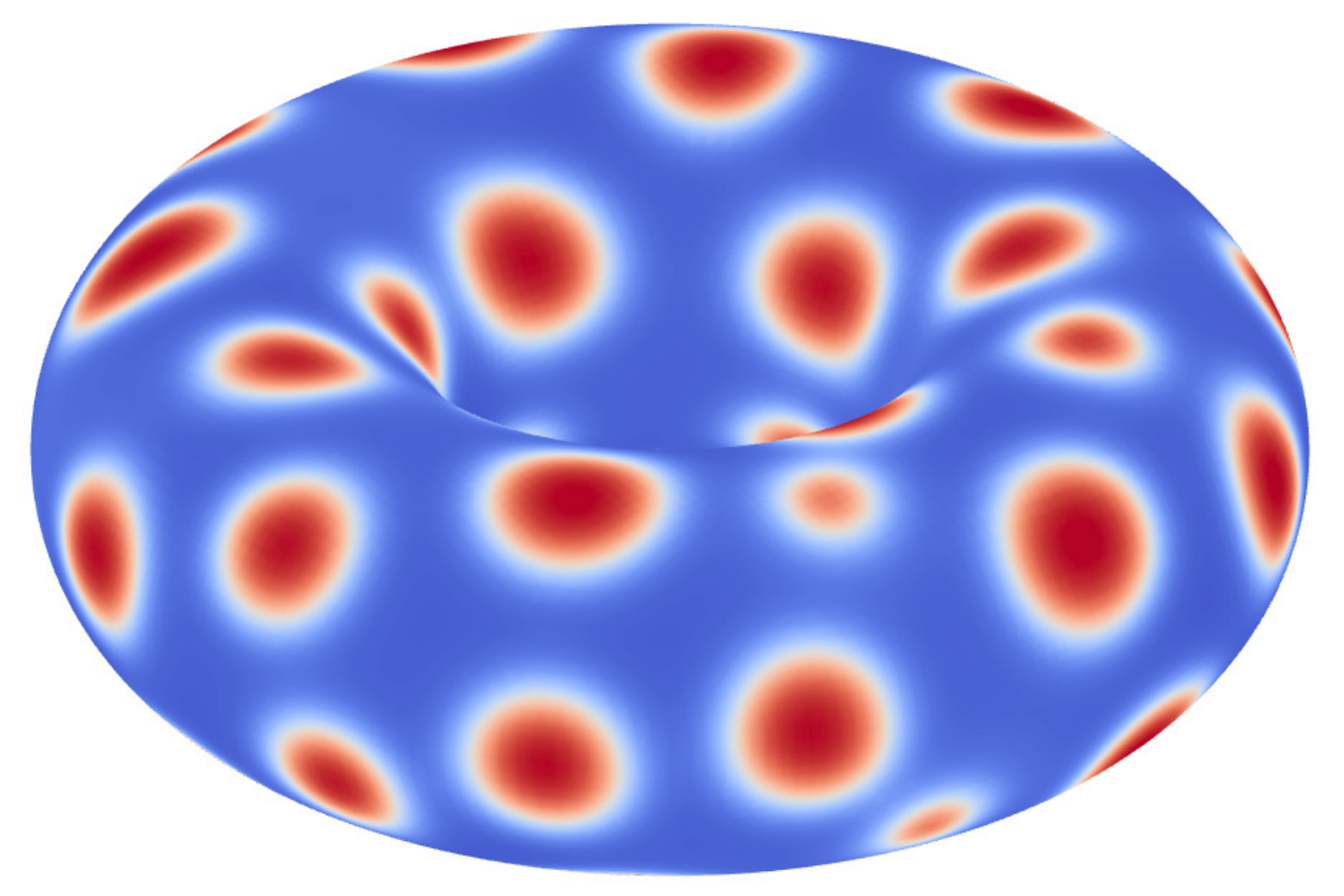} & \includegraphics[width=2.75cm]{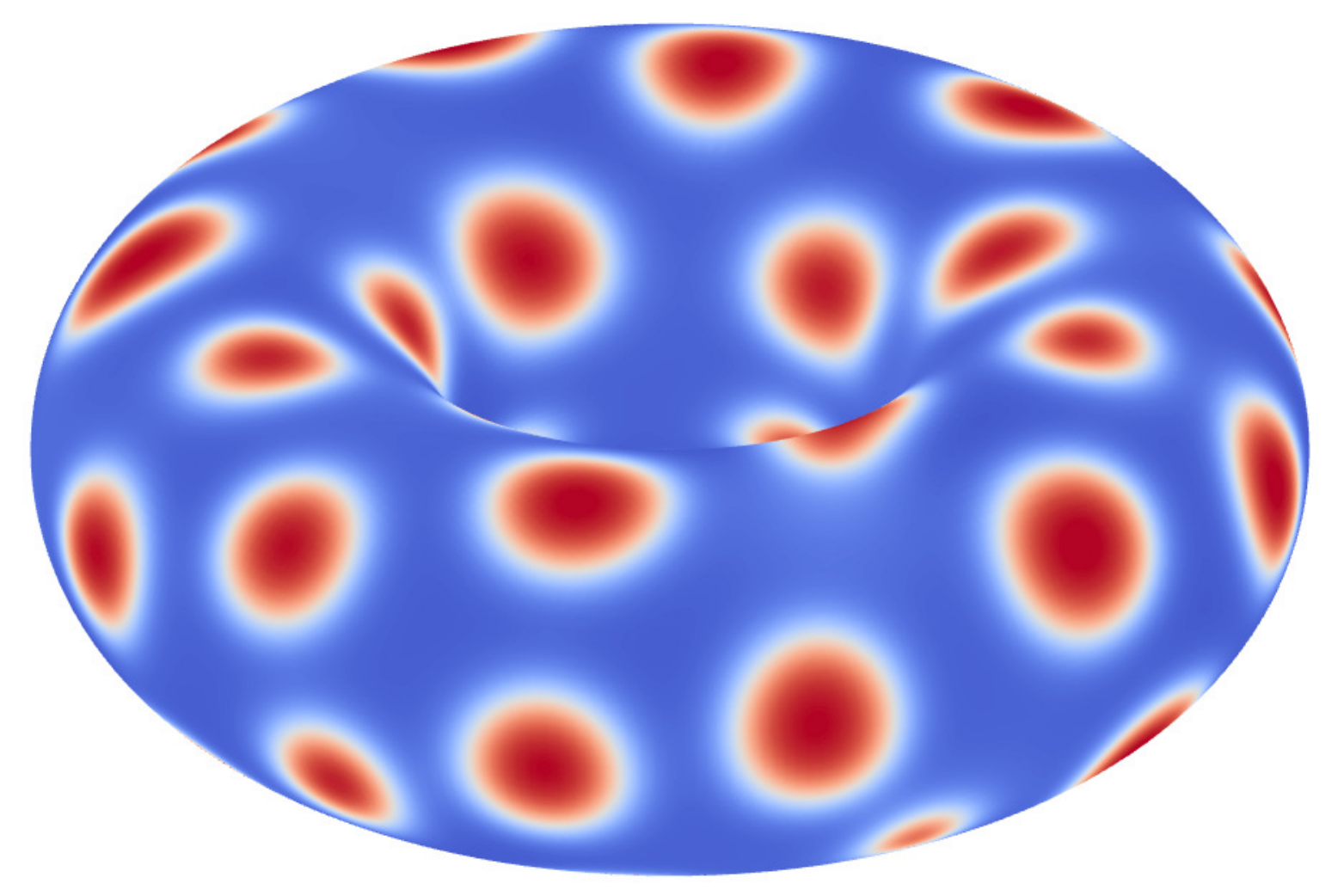} & \includegraphics[width=2.75cm]{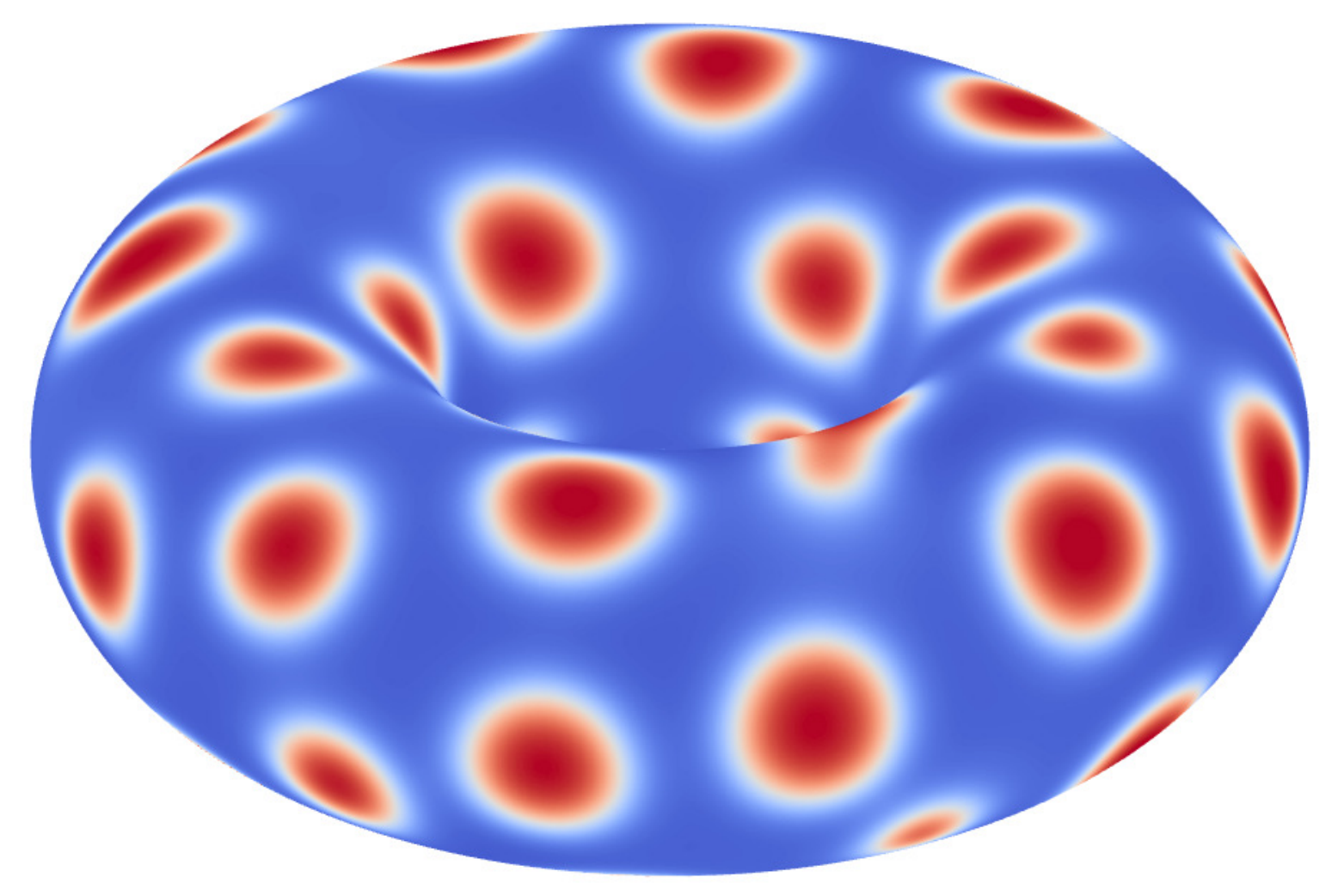} & \includegraphics[width=2.75cm]{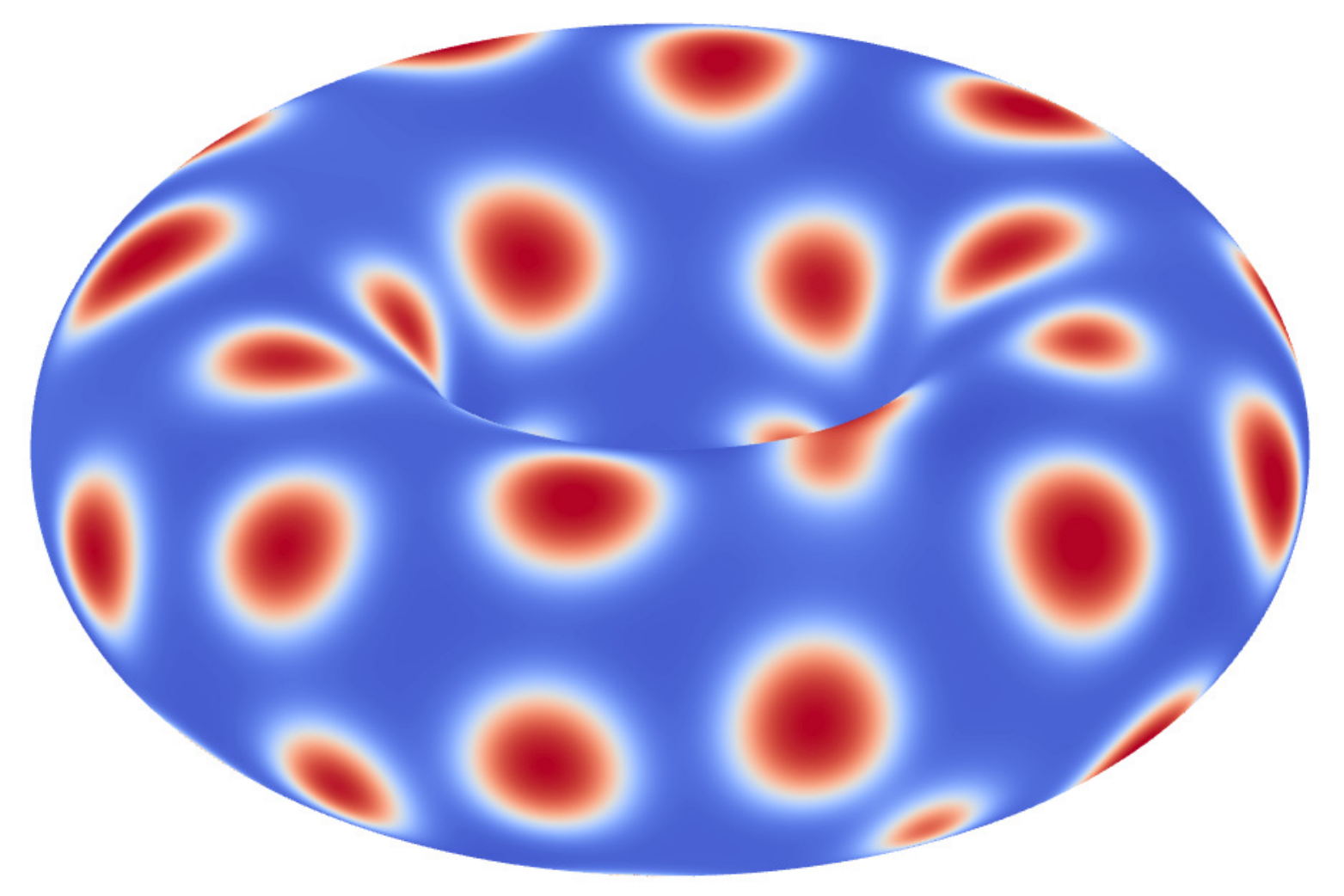} \tabularnewline
		1.0  & \includegraphics[width=2.75cm]{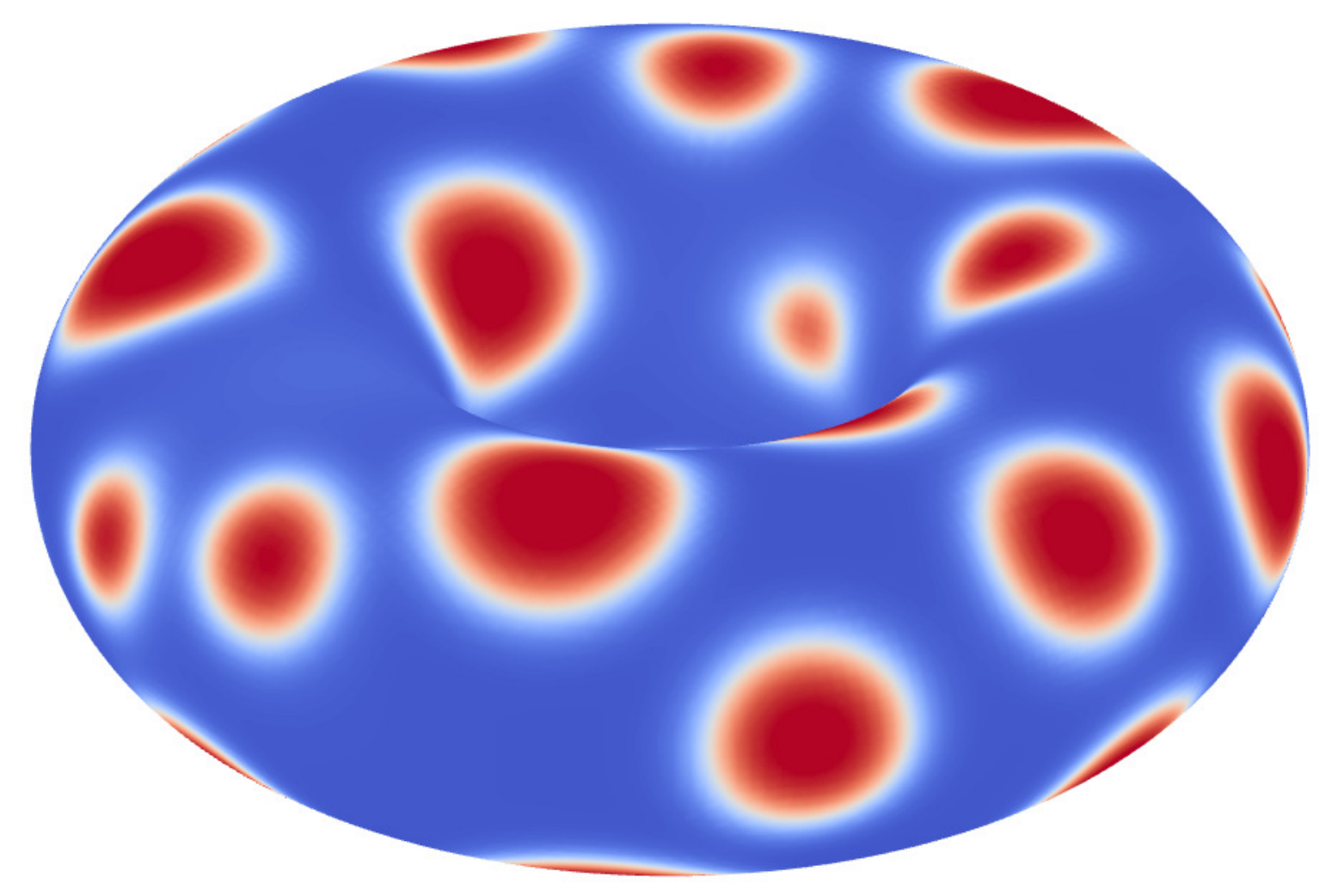} & \includegraphics[width=2.75cm]{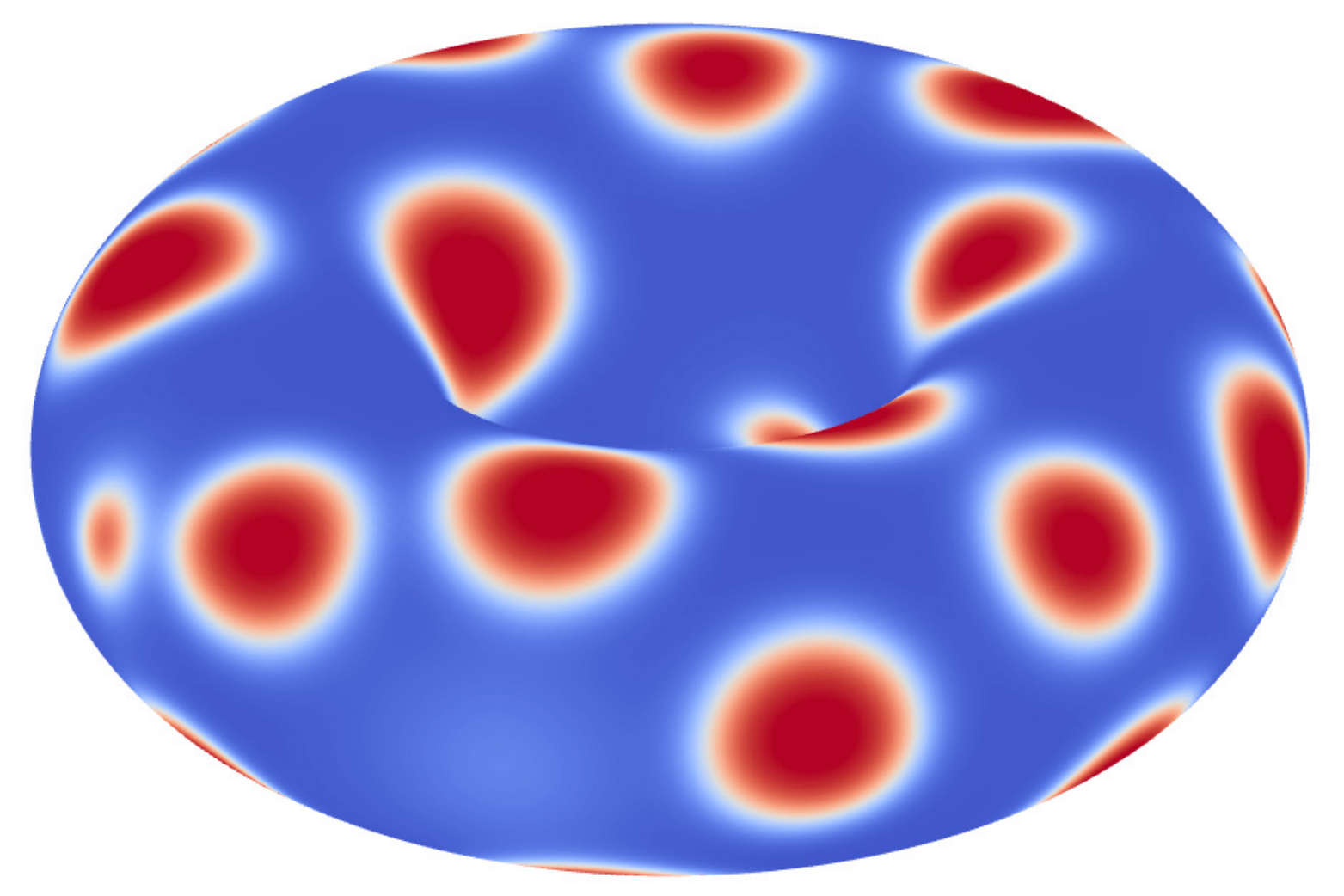} & \includegraphics[width=2.75cm]{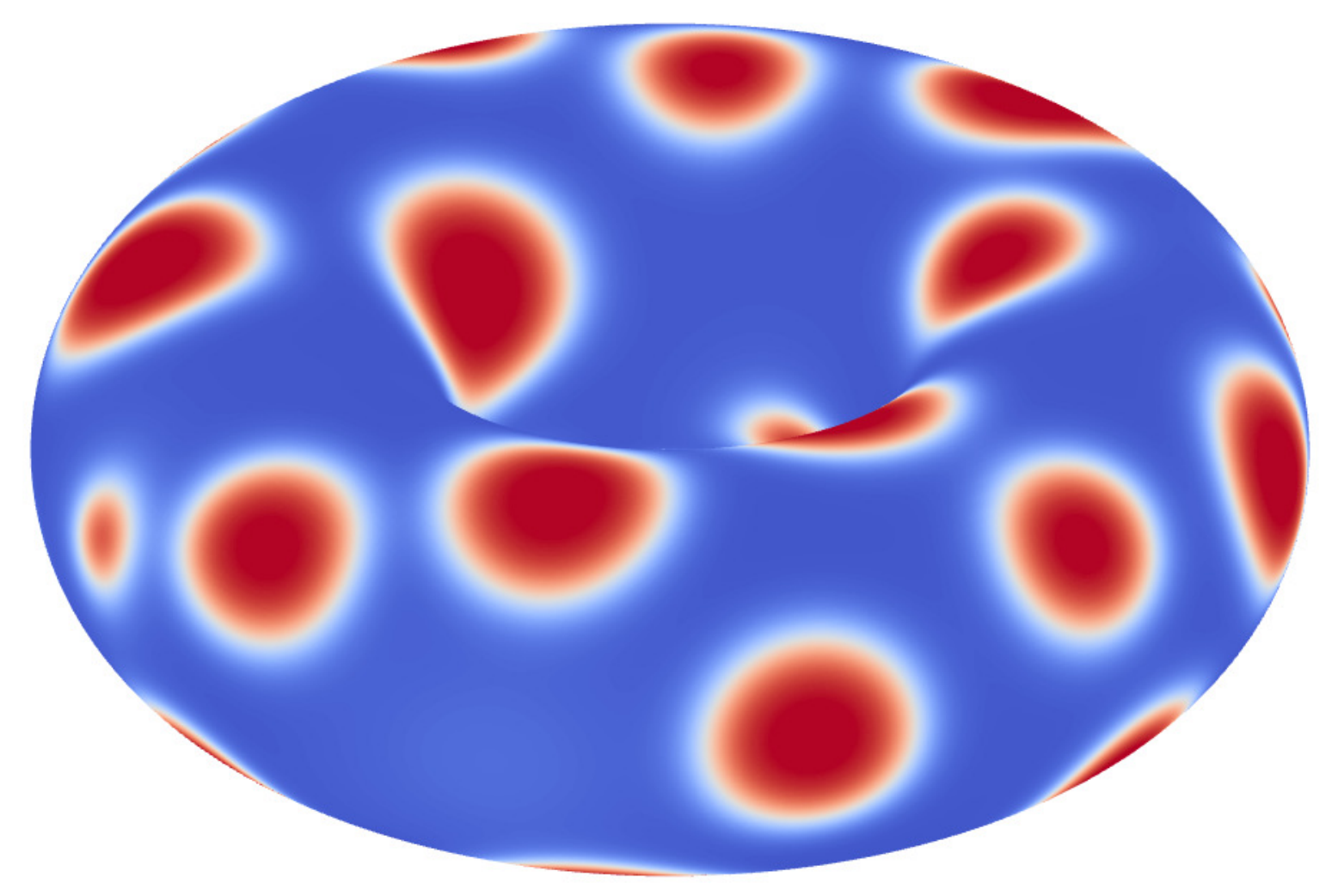} & \includegraphics[width=2.75cm]{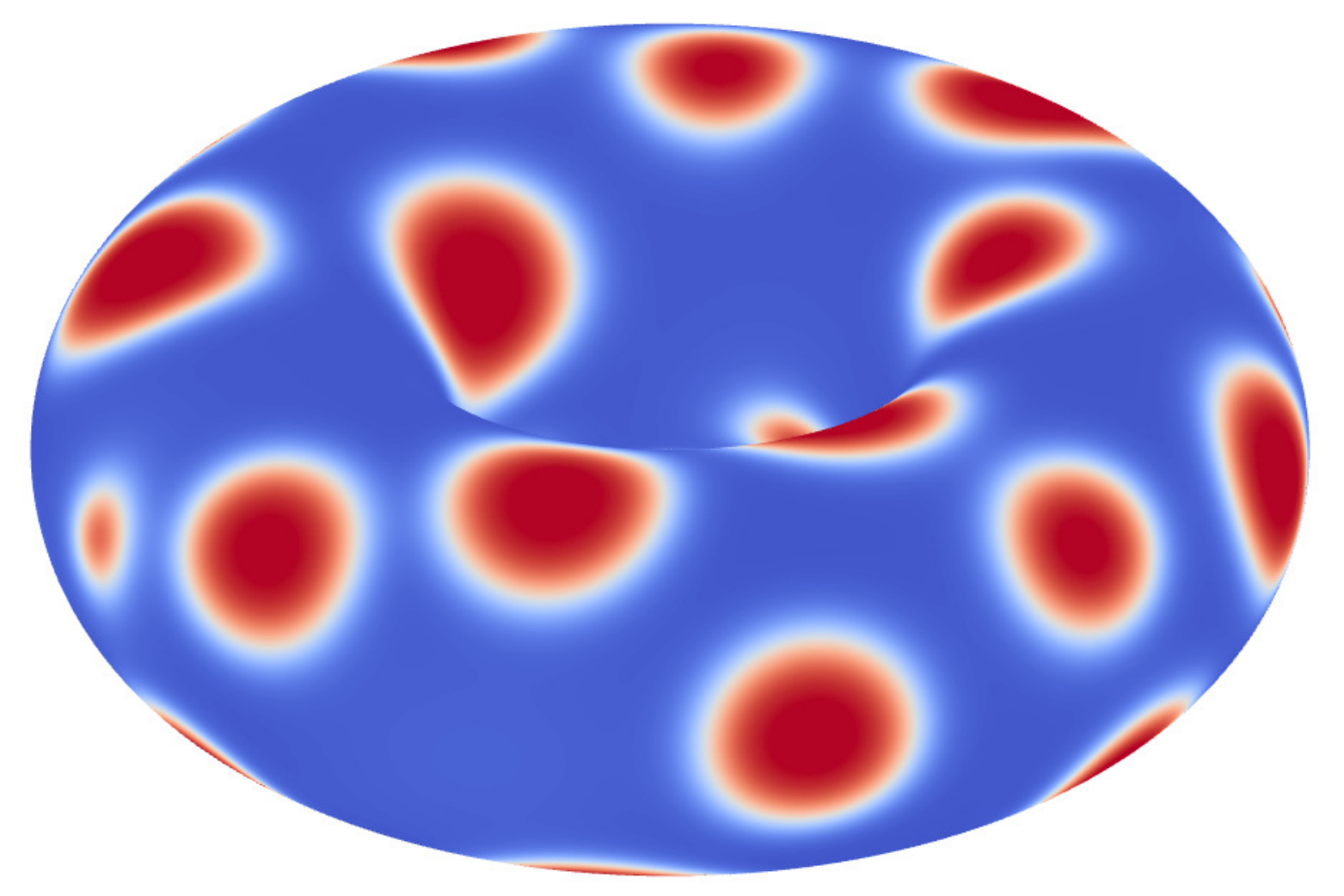} \tabularnewline
		2.0  & \includegraphics[width=2.75cm]{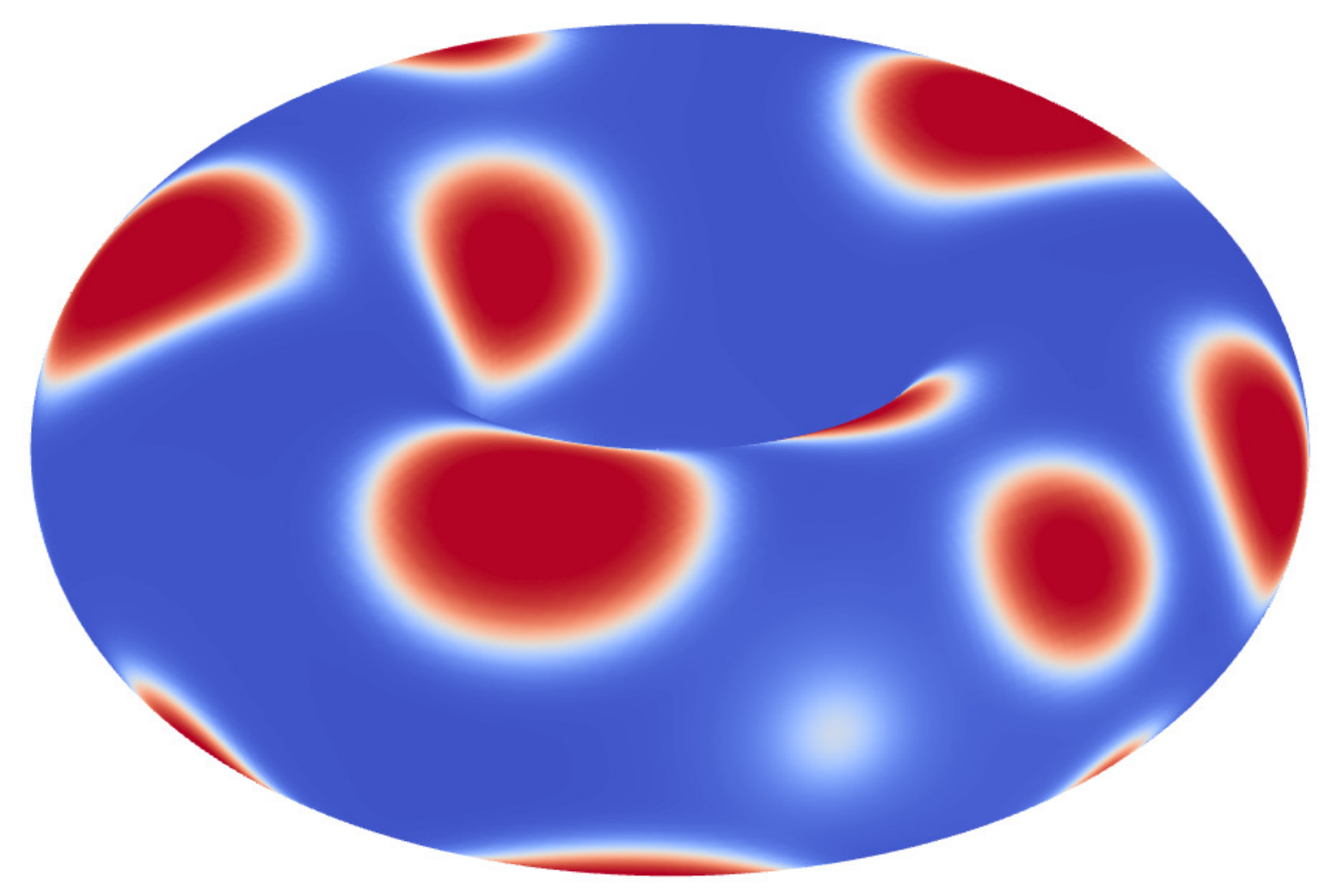} & \includegraphics[width=2.75cm]{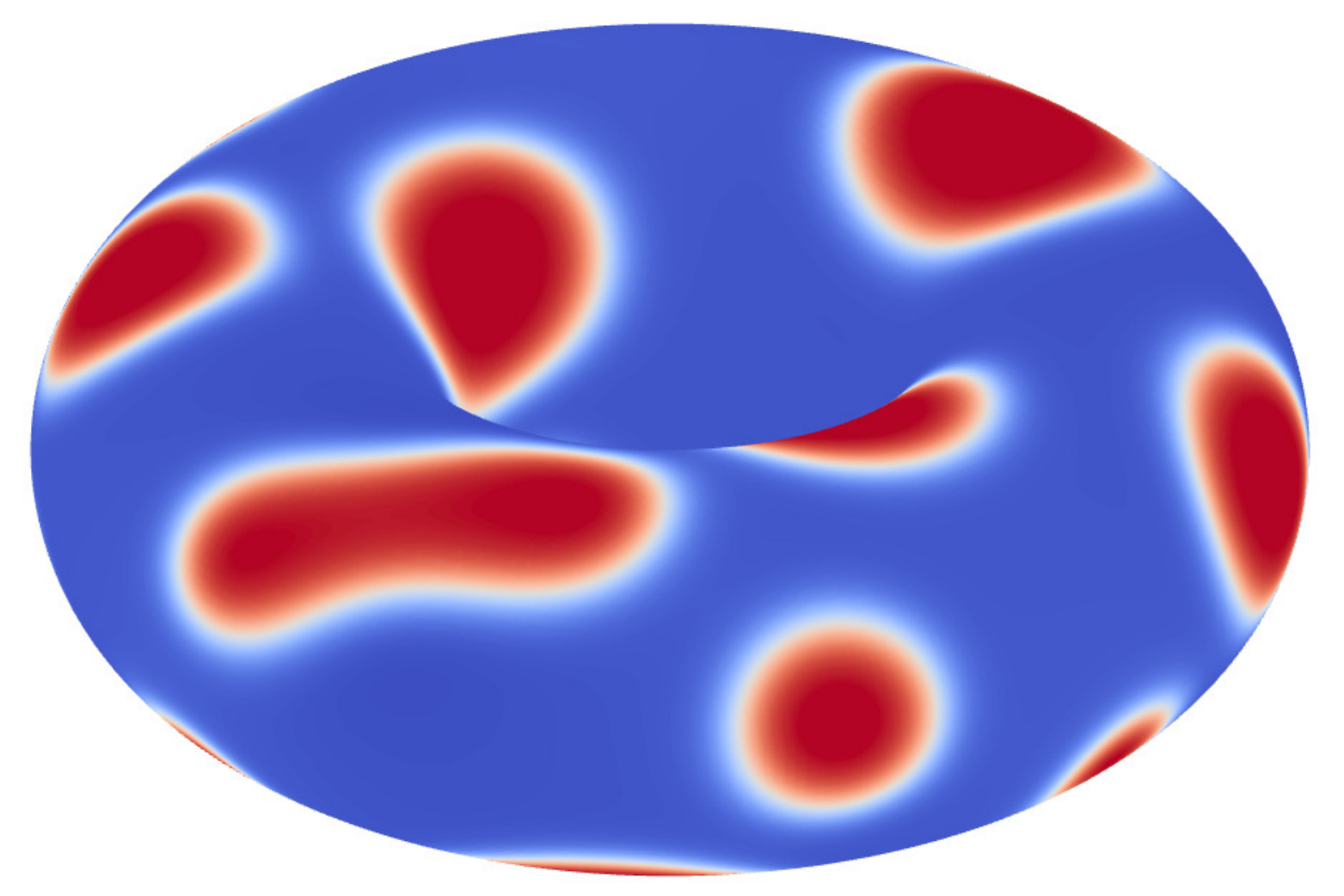} & \includegraphics[width=2.75cm]{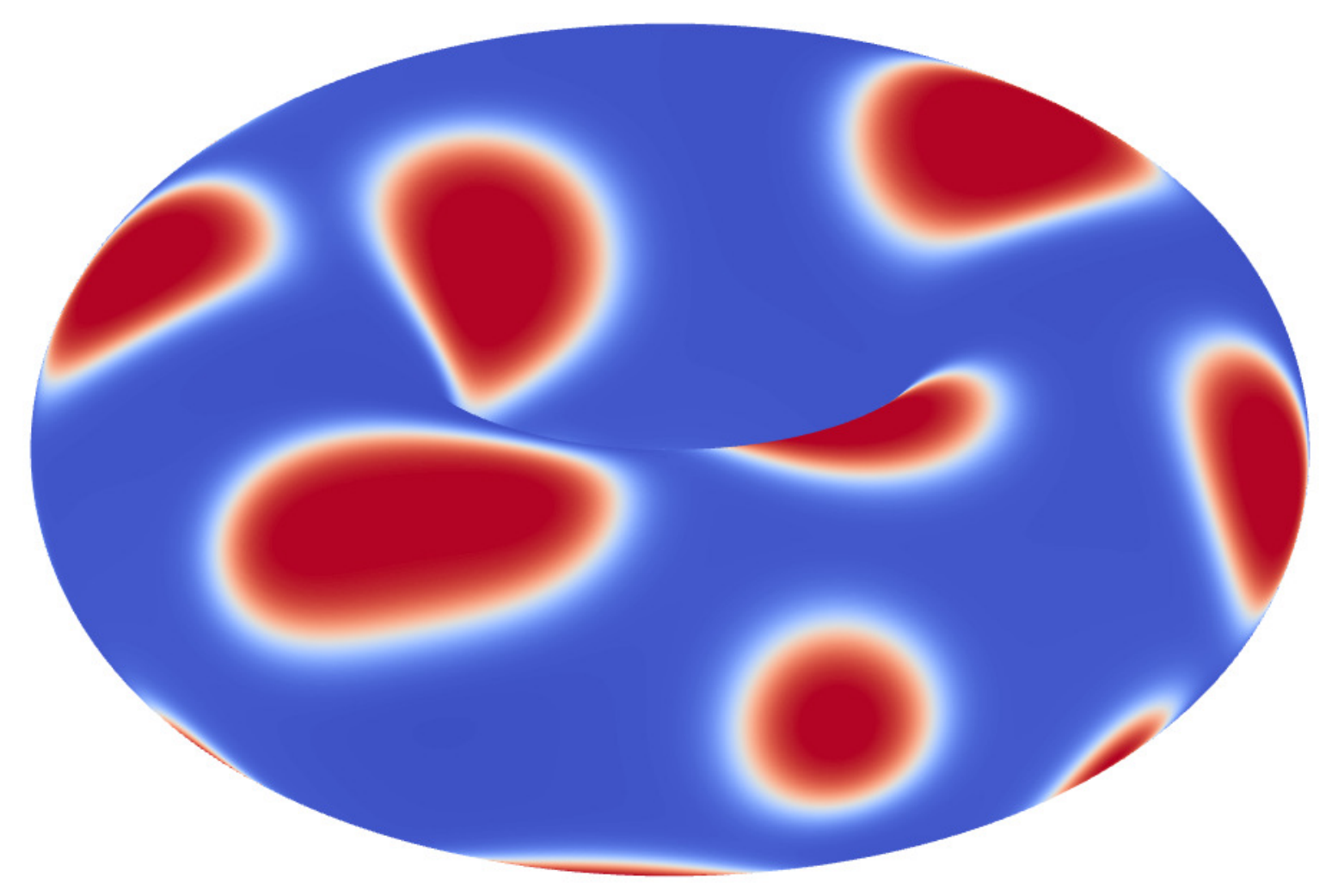} & \includegraphics[width=2.75cm]{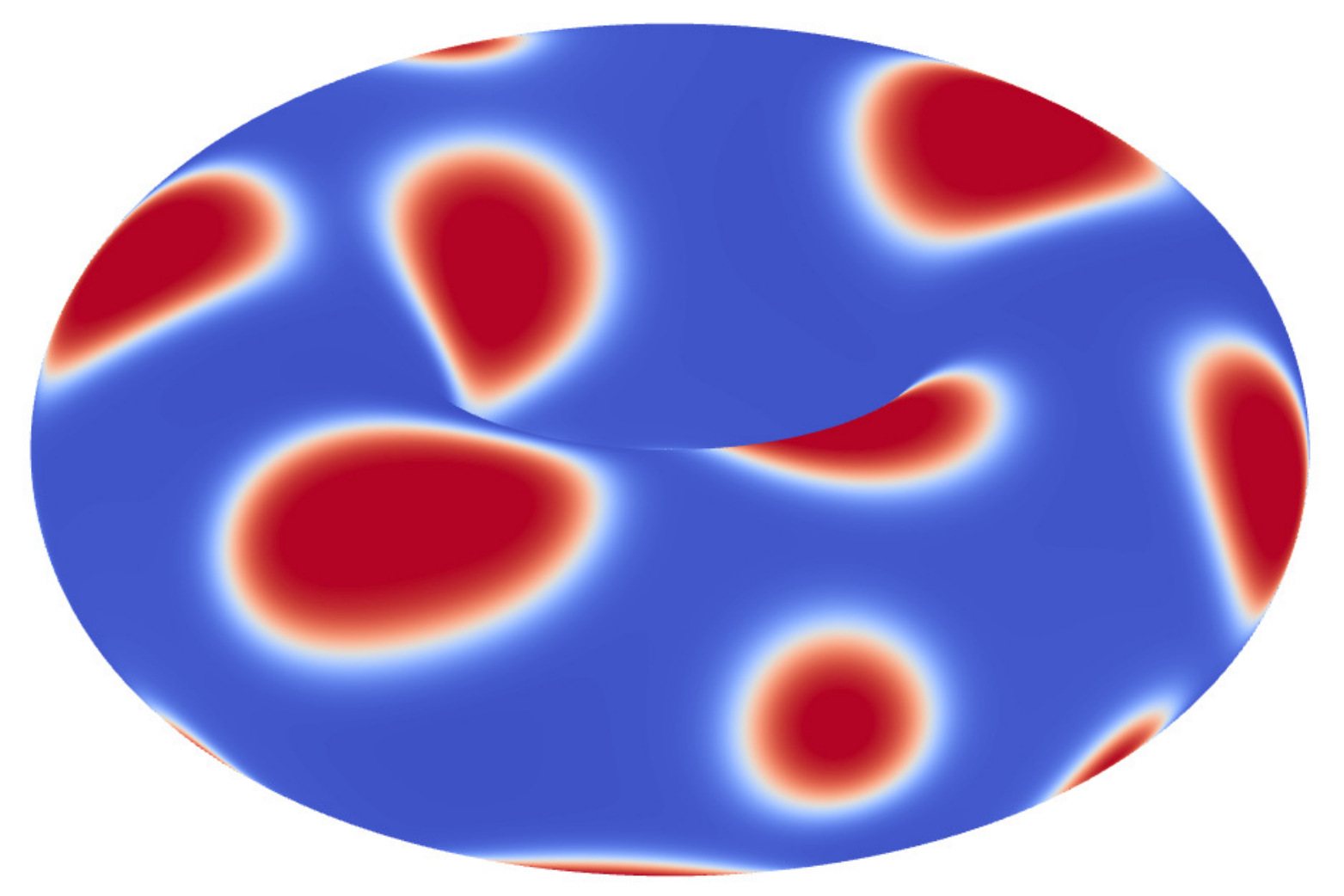} \tabularnewline
		5.0  & \includegraphics[width=2.75cm]{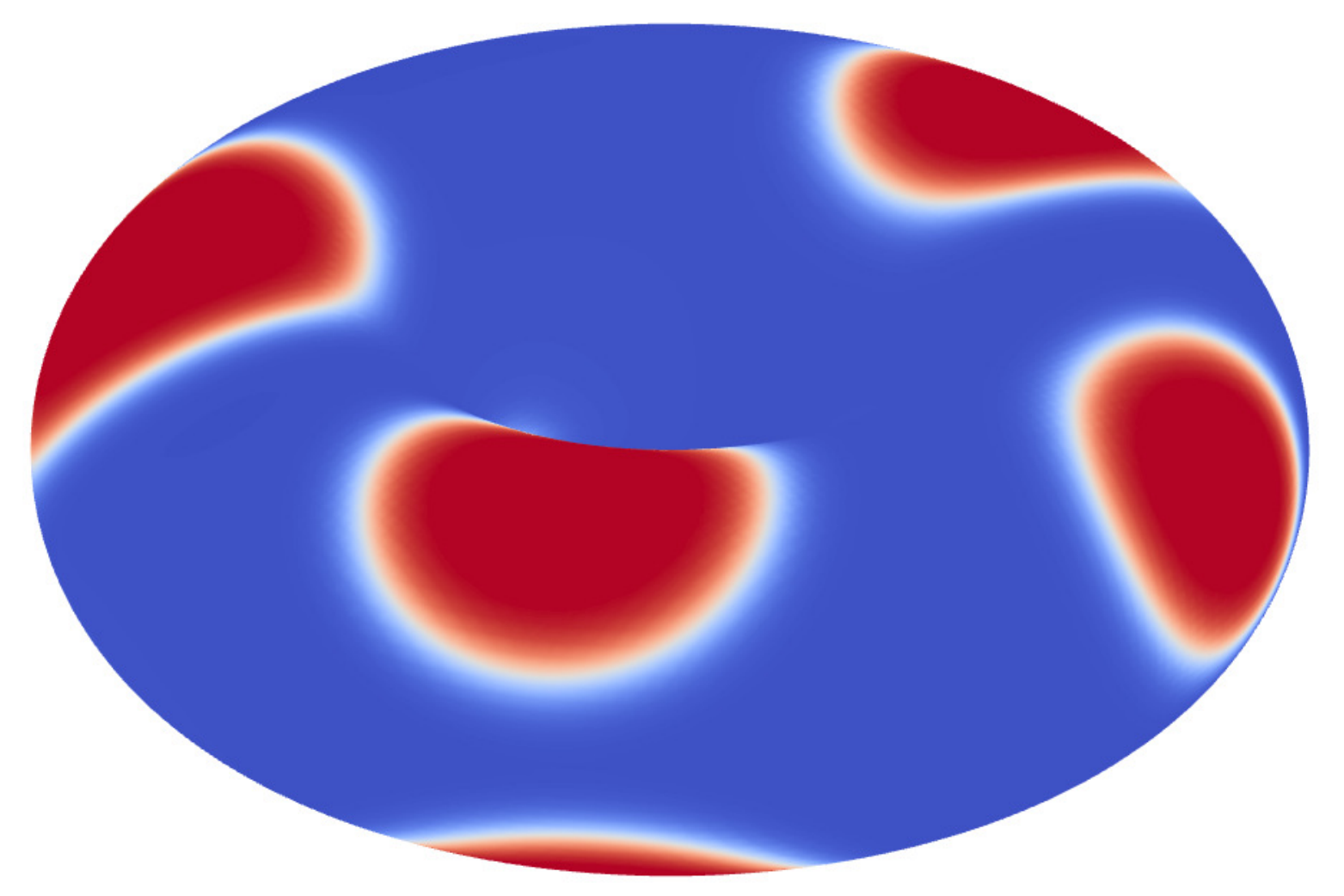} & \includegraphics[width=2.75cm]{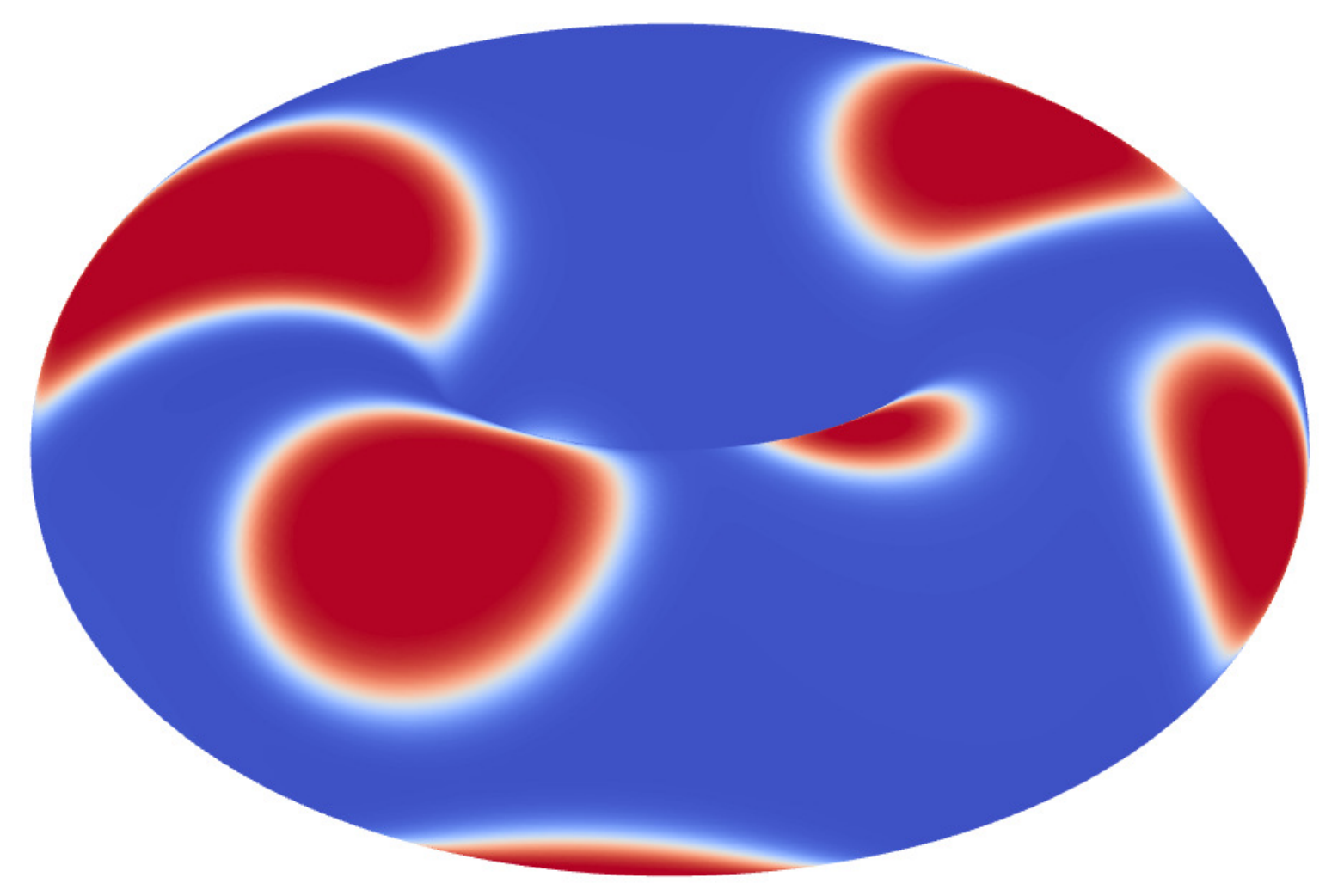} & \includegraphics[width=2.75cm]{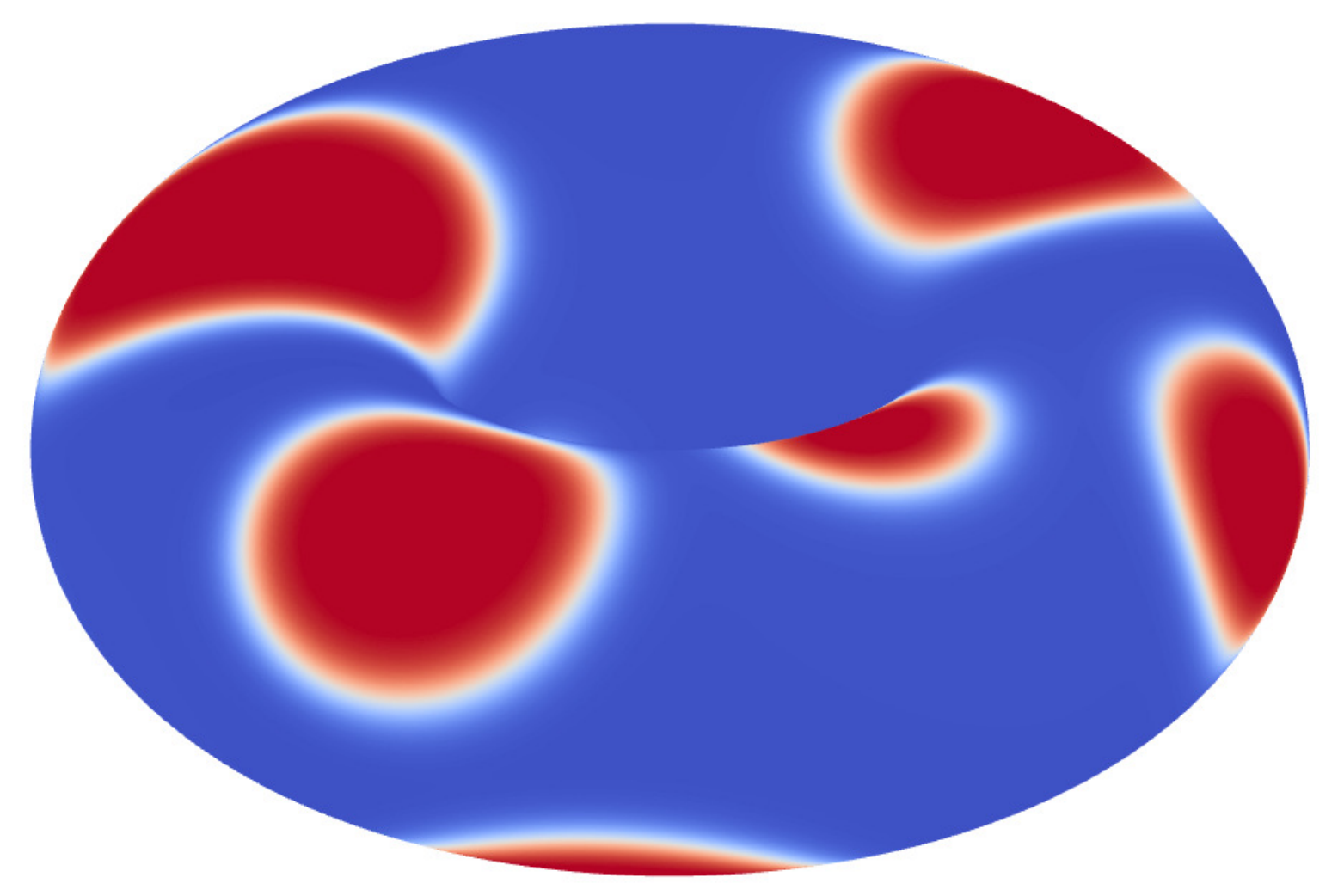} & \includegraphics[width=2.75cm]{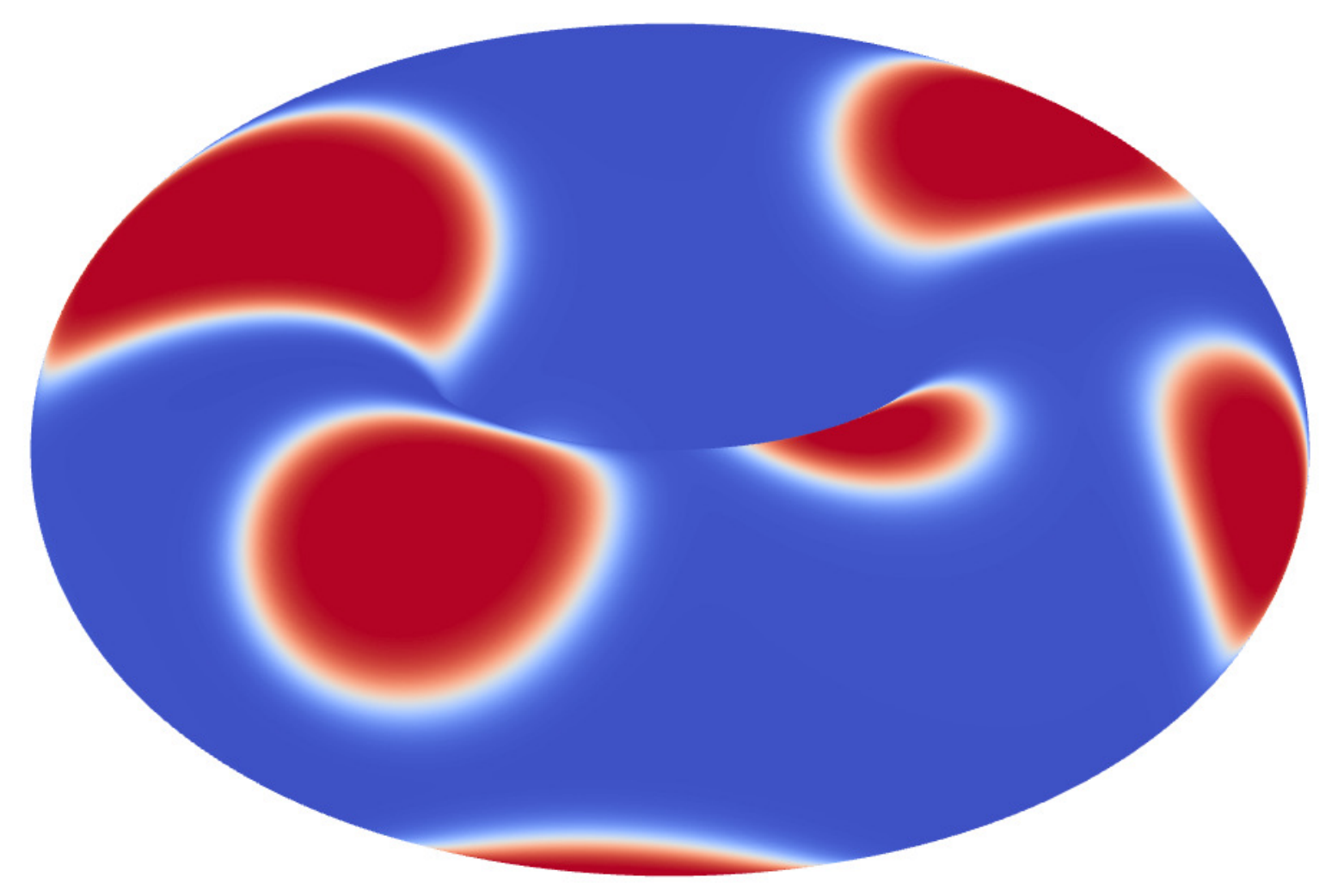} \tabularnewline
	\end{tabular}
	\caption{Evolution of the Cahn-Hilliard equation on a torus using a time step
	of $10^{-4}$ and varying grid sizes. The results do not qualitatively change when using a mesh size of $161^3$ or higher.}
	\label{fig:qualresults}
\end{figure}

The errors for both the BDF1 and BDF2 schemes as a function of time step for various grid spacings are shown in 
Fig. \ref{fig:timeConv}.
Several items become apparent from this figure. First, at 
large time steps both the first-order and second-order schemes converge at the expected
rate. Second, it is observed that for a given
grid size at sufficiently small time steps the error becomes constant. 
In fact, both the first- and second-order in time schemes have the same error
at small time steps. 

This behavior can be explained by considering the expected error. Due to the 
interpolation operators this solver is similar to those used for semi-Lagrangian
methods. In semi-Lagrangian methods the error obeys
$\epsilon\approx \mathcal{O}(\Delta t^k)+ \mathcal{O}(h^l)+ \mathcal{O}(h^p/\Delta t)$,
where $k$ is the order of
time discretization, $l$ is the truncation error of the derivative approximations, and 
$p$ is the truncation error of the interpolant~\cite{Velmurugan}. 
As the time-discretization errors decrease for a fixed grid spacing the overall
error should be dominated by the interpolation error, $\mathcal{O}(h^p/\Delta t)$.
It is interesting to note that the error for this particular solution does 
not increase with decreasing time step, but remains constant. 
One possible explanation is that the closest point location does not change
over the course of the current simulation, which could mean that the $\mathcal{O}(h^p/\Delta t)$ 
term is not present in the current error. Another possible explanation is that the time step 
simply has not become small enough for the effect to become evident.

The errors using $\Delta t=2.5\times 10^{-6}$ with the BDF2 scheme 
are shown in Table \ref{table:gridConv}.
At this time step the error is dominated by the grid spacing. 
It is unknown why the convergence rate between the $129^3$ and $161^3$ grids
is much lower than the other rates. Despite this, the spatial convergence is 
approximately second-order for both the phase field and chemical potential.

Next consider the qualitative convergence of the Cahn-Hilliard system. 
The Cahn-Hilliard system is solved on a torus with an inner radius
of 0.5 and an outer radius of 1.5 on a grid which spans $[-1.75,1.75]^3$. The mixing 
energy in this case is $g(f)=(f^4)/40 - (f^2)/20$.
To ensure that the initial condition is the same across all grid sizes,
a random initial concentration which spans $[-0.41,-0.39]$ is created on the $97^3$ grid and is then
interpolated using a spline interpolant
to the finer meshes. In this case the time step is fixed at $10^{-4}$
and the BDF2 scheme is used. 
The other parameters in this study are $\Cn=0.05$ and $\Pe=1.0$.
To ensure the conservation of the phase field over long simulation times
the correction of Xu, Li, Lowengrub, and Zhao has been implemented~\cite{xu2006level}.

The results for various times can be seen in Fig.~\ref{fig:qualresults}.
The initial homogeneous phase quickly
segregates, as seen at t=$0.2$, into many separate circular domains. After segregation the domains 
start to slowly coarsen in time. The large domains grow larger in size at the
expense of the smaller domains. From Fig.~\ref{fig:qualresults}, it can be
seen that the solution on the grid size of $97^3$ differs from the others.
For example, at a time of $t=2$, a domain which exists in the larger 
grids (in the lower-right corner of the figure) does not exist in the $97^3$ case. Additionally,
only five visible domains exist at a time of $t=5$ for the $97^3$ case, while six are seen in the 
other mesh sizes. From these results it appears that qualitative convergence is obtained
for meshes larger than $161^3$. This indicates that sufficient spatial resolution
is required to capture all aspects of the dynamics.

\section{Conclusions}

In this work the surface Cahn-Hilliard system is solved via a set of coupled second-order differential equations. The surface operators are discretized using
the Closest-Point Method. To aid in the solution of the resulting linear
systems, a Schur decomposition is used as a preconditioning matrix.
For large time steps the method converges at a rate dependent on the time discretization method. At small time
steps the error of the scheme becomes dominated by the spatial errors. Based on this result the second-order scheme was used to demonstrate
qualitative convergence of the Cahn-Hilliard system on a torus with a random initial condition.

\section*{References}

\bibliography{CHsolver}

\end{document}